\documentclass[twoside,11pt]{article}

\usepackage{amsmath,amssymb, amsthm}
\usepackage{graphicx}

\newtheorem{thm}{Theorem}
\newtheorem{cor}{Corollary}
\newtheorem{lem}{Lemma}

\theoremstyle{definition}

\theoremstyle{remark}
\newtheorem{rem}{Remark}
\numberwithin{equation}{section}
\begin{document}
\pagestyle{myheadings} \markboth{ \rm \centerline {R. N. Mohapatra and B. Szal}} {\rm \centerline { On trigonometric approximation of functions in the $L^{p}$ norm}}

\begin{titlepage}
\title{\bf {On trigonometric approximation of functions in the $L^{p}$ norm}}
\author {\bf\Large R. N. Mohapatra $^{1}$ and 
B. Szal $^{2}$\\
{\small $^1$ University of Central Florida}\\
{\small Orlando, FL 32816, USA}\\
{\small ramm@pegasus.cc.ucf.edu}\\
{\small $^2$ University of Zielona G\'{o}ra}\\
{\small Faculty of Mathematics, Computer Science and Econometrics}\\
{\small 65-516 Zielona G\'{o}ra, ul. Szafrana 4a, Poland}\\
{\small B.Szal@wmie.uz.zgora.pl} }
\end{titlepage}

\date{}
\maketitle
\begin{abstract}In this paper we obtain degree of approximation of functions in $L^{p}$ by
operators associated with their Fourier series using integral modulus of
continuity. These results generalize many know results and are proved under
less stringent conditions on the infinite matrix.
\end{abstract}

\noindent{\it Keywords and phrases:} Class $Lip\left( \alpha ,p\right) $; Trigonometric
approximation; $L^{p}-$norm.

\noindent { \it 2000 Mathematics Subject Classification:}  42A10, 41A25.

\maketitle

\section{Introduction}

Let $f$ be $2\pi $ periodic and $f\in L^{p}\left[ 0,2\pi \right] $ for $%
p\geq 1$. Denote by 
\begin{equation*}
S_{n}\left( f\right) =S_{n}\left( f;x\right) =\frac{a_{0}}{2}%
+\sum\limits_{k=1}^{n}\left( a_{k}\cos kx+b_{k}\sin kx\right) \equiv
\sum\limits_{k=0}^{n}U_{k}\left( f;x\right) 
\end{equation*}%
partial sum of the first $\left( n+1\right) $ terms of the Fourier series of 
$f\in L^{p}$ $\left( p\geq 1\right) $ at a point $x$, and by%
\begin{equation*}
\omega _{p}\left( f;\delta \right) =\underset{0<\left\vert h\right\vert \leq
\delta }{\sup }\left\{ \frac{1}{2\pi }\int\limits_{0}^{2\pi }\left\vert
f\left( x+h\right) -f\left( x\right) \right\vert ^{p}dx\right\} ^{\frac{1}{p}%
}
\end{equation*}%
the integral modulus of continuity of $f\in L^{p}$. If, for $\alpha >0$, $%
\omega _{p}\left( f;\delta \right) =O\left( \delta ^{\alpha }\right) $, then
we write $f\in Lip\left( \alpha ,p\right) $ $\left( p\geq 1\right) $.

Throughout $\left\Vert \cdot \right\Vert _{L^{p}}$ will denote $L^{p}-$norm,
defined by%
\begin{equation*}
\left\Vert f\right\Vert _{L^{p}}=\left\{ \frac{1}{2\pi }\int\limits_{0}^{2%
\pi }\left\vert f\left( x\right) \right\vert ^{p}dx\right\} ^{\frac{1}{p}}%
\text{ \ \ }\left( f\in L^{p}\left( p\geq 1\right) \right) .
\end{equation*}%
In the present paper, we shall consider approximation of $f\in L^{p}$ by
trigonometrical polynomials $T_{n}\left( f;x\right) $, where%
\begin{equation*}
T_{n}\left( f;x\right) =T_{n}\left( f,A;x\right)
:=\sum_{k=0}^{n}a_{n,k}S_{k}\left( f;x\right) \text{ \ \ \ }\left(
n=0,1,2,...\right) 
\end{equation*}%
and $A:=\left( a_{n,k}\right) $ \ be a lower triangular infinite matrix of
real numbers such that:

\begin{equation}
a_{n,k}\geq 0\text{ for }k\leq n\text{ and }a_{n,k}=0\text{ for }k>n\text{ }%
\left( k,n=0,1,2....\right)   \label{a1}
\end{equation}%
and%
\begin{equation}
\sum_{k=0}^{n}a_{n,k}=1\left( n=0,1,2...\right) .  \label{a2}
\end{equation}

If $a_{n,k}=\frac{p_{k}}{P_{m}}$, where $P_{n}=p_{0}+p_{1}+...+p_{n}\neq 0$ $%
\left( n\geq 0\right) $, then we shall call this trigonometrical polynomials
by%
\begin{equation*}
R_{n}\left( f;x\right) =\frac{1}{P_{n}}\sum\limits_{k=0}^{n}p_{k}S_{k}%
\left( f;x\right) \text{ \ \ \ \ }\left( n=0,1,2....\right) .
\end{equation*}

The case $a_{n,k}=\frac{1}{n+1}$ for $k\leq n$ and $a_{n,k}=0$ for $k>n$ of $%
T_{n}\left( f;x\right) $ yields%
\begin{equation*}
\sigma _{n}\left( f;x\right) =\frac{1}{n+1}\sum\limits_{k=0}^{n}S_{k}\left(
f;x\right) \text{ \ \ \ \ }\left( n=0,1,2....\right) .
\end{equation*}%
We shall also use the notations%
\begin{equation*}
\Delta a_{k}=a_{k}-a_{k+1}\text{, \ \ \ \ \ }\Delta
_{k}a_{n,k}=a_{n,k}-a_{n,k+1}
\end{equation*}%
and we shall write $I_{1}\ll I_{2}$ if there exists a positive constant $K$
such that $I_{1}\leq KI_{2}$.

Let $C:=\left( C_{n}\right) =\frac{1}{n+1}\sum\limits_{k=0}^{n}c_{k}$,
where $c:=\left( c_{n}\right) $ is a sequence of nonnegative numbers. The
sequence $c$ is called a nondecreasing (nonincreasing) mean sequence, briefly $%
NDMS$ ($NIMS$), if $C\in NDS$ ($C\in NIS$), where $NDS$ ($NIS$) is the class
of nonnegative and nondecreasing (nonincreasing) sequences.

A nonnegative sequence $c:=\left( c_{n}\right) $ is called almost monotone
decreasing ($AMDS$) (increasing ($AMIS$)) if there exists a constant $%
K:=K\left( c\right) $, depending on the sequence $c$ only, such that for all 
$n\geq m$%
\begin{equation*}
c_{n}\leq Kc_{m}\text{ \ \ \ \ }\left( Kc_{n}\geq c_{m}\right) .
\end{equation*}%
Such sequences will be denoted by $c\in AMDS$ and $c\in AMIS$, respectively.

If $C\in AMDS$ ($C\in AMIS)$, then we shall say that $c$ is almost monotone
decreasing (increasing) mean sequence, briefly $c\in AMDMS$ $\left( c\in
AMIMS\right) $.

When we write that a sequence $\left( a_{n,k}\right) $ belongs to one of the
above classes, it means that it satisfies the required conditions from the above definitions with
respect to $k=0,1,2,...,n$ for all $n.$

A sequence $c:=\left( c_{n}\right) $ of nonnegative numbers tending to zero
is called a rest bounded variation sequence (rest bounded variation mean
sequence), or briefly $c\in RBVS$ ($c\in RBVMS$), if it has the property%
\begin{equation}
\sum\limits_{k=m}^{\infty }\left\vert \Delta c_{k}\right\vert \leq K\left(
c\right) c_{m}\text{ \ \ \ \ }\left( \sum\limits_{k=m}^{\infty }\left\vert
\Delta C_{k}\right\vert \leq K\left( c\right) C_{m}\right)  \label{1}
\end{equation}%
for all natural numbers $m$, where $K\left( c\right) $ is a constant
depending only on $c$.

A sequence $c:=\left( c_{n}\right) $ of nonnegative numbers will be called
a head bounded variation sequence (head bounded variation mean sequence),
or briefly $c\in HBVS$ ($c\in HBVMS$), if it has the property%
\begin{equation}
\sum\limits_{k=0}^{m-1}\left\vert \Delta c_{k}\right\vert \leq K\left(
c\right) c_{m}\text{ \ \ \ \ }\left( \sum\limits_{k=0}^{m-1}\left\vert
\Delta C_{k}\right\vert \leq K\left( c\right) C_{m}\right)   \label{2}
\end{equation}%
for all natural numbers $m$, or only for all $m\leq N$ if the sequence $c$
has only a finite number of nonzero terms and the last nonzero terms is $c_{N}$.

Therefore we assume that the sequence $\left( K\left( \alpha _{n}\right)
\right) _{n=0}^{\infty }$ is bounded, that is, there exists a constant $%
K$ such that%
\begin{equation*}
0\leq K\left( \alpha _{n}\right) \leq K
\end{equation*}%
holds for all $n$, where $K\left( \alpha _{n}\right) $ denote the sequence
of constants appearing in the inequalities (\ref{1}) or (\ref{2}) for the
sequence $\alpha _{n}:=\left( a_{nk}\right) _{k=0}^{\infty }$.Now we can
mention the conditions to be used later on. Let $A_{n,m}=\frac{1}{m+1}%
\sum\limits_{k=0}^{m}a_{n,k}$. We assume that for all $n$ and $0\leq m\leq
n $%
\begin{equation*}
\sum\limits_{k=m}^{\infty }\left\vert \Delta _{k}a_{nk}\right\vert \leq
Ka_{nm}\text{ \ \ \ \ }\left( \sum\limits_{k=m}^{\infty }\left\vert \Delta
_{k}A_{nk}\right\vert \leq KA_{nm}\right)
\end{equation*}%
and%
\begin{equation*}
\sum\limits_{k=0}^{m-1}\left\vert \Delta _{k}a_{nk}\right\vert \leq Ka_{nm}%
\text{ \ \ \ \ }\left( \sum\limits_{k=0}^{m-1}\left\vert \Delta
_{k}A_{nk}\right\vert \leq KA_{nm}\right)
\end{equation*}%
hold if $\alpha _{n}:=\left( a_{nk}\right) _{k=0}^{\infty }$ belongs to $%
RBVS $ ($RBVMS$) or $HBVS$ ($HBVMS$), respectively.

It is clear that%
\begin{eqnarray*}
NIS &\subset &RBVS\subset AMDS, \\
NIMS &\subset &RBVMS\subset AMDMS
\end{eqnarray*}%
and%
\begin{eqnarray*}
NDS &\subset &HBVS\subset AMIS, \\
NDMS &\subset &HBVMS\subset AMIMS.
\end{eqnarray*}%
In the present paper we shall show that $NIS\subset NIMS$, $AMDS\subset AMDMS
$, $NDS\subset NDMS$ and $AMIS\subset AMIMS$, too.

In 1937 E. Quade \cite{3} proved that, if $f\in Lip\left( \alpha ,p\right) $
for $0<\alpha \leq 1$, then $\left\Vert \sigma _{n}\left( f\right)
-f\right\Vert _{L^{p}}=O\left( n^{-a}\right) $ for either $p>1$ and $%
0<\alpha \leq 1$ or $p=1$ and $0<\alpha <1$. He also showed that, if $%
p=\alpha =1$, then $\left\Vert \sigma _{n}\left( f\right) -f\right\Vert
_{L^{1}}=O\left( n^{-1}\log \left( n+1\right) \right) $.

There are several generalizations of the above result for $p>1$ (see, for
example\cite{4, 5, 6}, \cite{7} and \cite{8}). In \cite{1} P. Chandra
extended the work of E. Quade and proved the following theorems:

\begin{thm}
Let $f\in Lip\left( \alpha ,p\right) $ and let $\left( p_{n}\right) $ be
positive. Suppose that either

$\left( i\right) $ $p>1$, $0<\alpha \leq 1$, and

$\left( ii\right) $ $\sum\limits_{k=0}^{n-1}\left\vert \Delta \left( \frac{%
P_{k}}{k+1}\right) \right\vert =O\left( \frac{P_{n}}{n+1}\right) $, or

$\left( i\right) $ $p=1$, $0<\alpha <1$, and

$\left( ii\right) $ $\left( p_{n}\right) $ is nondecreasing and

\begin{equation}
\left( n+1\right) p_{n}=O\left( P_{n}\right) .  \label{3}
\end{equation}

Then%
\begin{equation*}
\left\Vert R_{n}\left( f\right) -f\right\Vert _{L^{p}}=O\left( n^{-\alpha
}\right) \text{.}
\end{equation*}
\end{thm}

\begin{thm}
Let $f\in Lip\left( 1,1\right) $ and let $\left( p_{n}\right) $ with (\ref{3}%
) be positive, and that%
\begin{equation*}
\left( \left( n+1\right) ^{\eta }p_{n}\right) \in NDS\text{ for some }\eta
>0.
\end{equation*}%
Then%
\begin{equation*}
\left\Vert R_{n}\left( f\right) -f\right\Vert _{L^{1}}=O\left( n^{-1}\right)
.
\end{equation*}
\end{thm}

In \cite{2} M. Mittal,B. Rhoades,V. Mishra and U. Singh obtained the same
degree of approximation as in above theorems, for a more general class of
lower triangular matrices, and deduced some of the results of P. Chandra.
Namely, they proved the following theorem:

\begin{thm}
Let $f\in Lip\left( \alpha ,p\right) $, and let $a_{n,k}\geq 0$ $\left(
k,n=0,1,...\right) $, $\left( a_{nk}\right) \in NDS$ or $\left(
a_{n,k}\right) \in NIS$ and%
\begin{equation*}
\left\vert \sum\limits_{k=0}^{n}a_{n,k}-1\right\vert =O\left( n^{-\alpha
}\right) .
\end{equation*}

$\left( i\right) $ If $p>1$, $0<\alpha <1$, $\left( n+1\right) \max \left\{
a_{n,0},a_{n,r}\right\} =O\left( 1\right) $, where $r:=\left[ \frac{n}{2}%
\right] $, then%
\begin{equation}
\left\Vert T_{n}\left( f\right) -f\right\Vert _{L^{p}}=O\left( n^{-\alpha
}\right) .  \label{4}
\end{equation}

$\left( ii\right) $ If $p>1$, $\alpha =1$, then (\ref{4}) is satisfied.

$\left( iii\right) $ If $p=1$, $0<\alpha <1$, and $\left( n+1\right) \max
\left\{ a_{n,0},a_{nn}\right\} =O\left( 1\right) $, then (\ref{4}) is
satisfied.
\end{thm}

In this paper we shall prove that the above mentioned theorems are valid with 
less stringent assumptions.

\section{Statement of the results}
\ Our first theroem deals with a number of embedding results.
\begin{thm}
The following embedding relations are valid:

$\left( i\right) $ $NIS\subset NIMS$,

$\left( ii\right) $ $NDS\subset NDMS,$

$\left( iii\right) $ $AMDS\subset AMDMS$,

$\left( iv\right) $ $AMIS\subset AMIMS.$
\end{thm}
\ Our next theorem deals with degree of convergence of operators involving the infinite matrix.
\begin{thm}
Let $f\in Lip\left( \alpha ,p\right) $ and (\ref{a1}), (\ref{a2}) hold. If
one of the conditions

$\left( i\right) $ $p>1$, $0<\alpha <1$ and $\left( a_{n,k}\right) \in AMIMS$%
,

$\left( ii\right) $ $p>1$, $0<\alpha <1$, $\left( a_{n,k}\right) \in AMDMS$
and $\left( n+1\right) a_{n,0}=O\left( 1\right) ,$

$\left( iii\right) $ $p>1$, $\alpha =1$ and $\sum\limits_{k=0}^{n-1}\left%
\vert \Delta _{k}A_{n,k}\right\vert =O\left( n^{-1}\right) $,

$\left( iv\right) $ $p=1$, $0<\alpha <1$, $\sum\limits_{k=0}^{n-1}\left%
\vert \Delta _{k}a_{n,k}\right\vert =O\left( n^{-1}\right) $ and $\left(
n+1\right) a_{n,n}=O\left( 1\right) ,$

$\left( v\right) $ $p=1$, $0<\alpha <1$, $\left( a_{n,k}\right) \in RBVS$
and $\left( n+1\right) a_{n,0}=O\left( 1\right) ,$

$\left( vi\right) $ $p=\alpha =1$, $\left( \left( k+1\right) ^{-\beta
}a_{n,k}\right) \in HBVS$ for some $\beta >0$ and $\left( n+1\right)
a_{n,n}=O\left( 1\right) $

maintains, then%
\begin{equation}
\left\Vert T_{n}\left( f\right) -f\right\Vert _{L^{p}}=O\left( n^{-\alpha
}\right) .  \label{5}
\end{equation}
\end{thm}

\begin{rem}
Let $f\in Lip\left( \alpha ,p\right) $, (\ref{a1}) and 
\begin{equation*}
\left\vert \sum\limits_{k=0}^{n}a_{n,k}-1\right\vert =O\left( n^{-\alpha
}\right) 
\end{equation*}%
hold. Under the assumptions of Theorem 5 $\left( i\right) -\left( vi\right) $
we can observe that the estimation (\ref{5}) are true, too.
\end{rem}

In the special cases, putting $a_{n,k}=\frac{p_{k}}{P_{n}}$, where $%
P_{n}=p_{0}+p_{1}+...+p_{n}\neq 0$ , we can derive from Theorem 5 the
following corollary:

\begin{cor}
Let $f\in Lip\left( \alpha ,p\right) $ and let $\left( p_{k}\right) $ be
positive. If one of the conditions

$\left( i\right) $ $p>1$, $0<\alpha <1$ and $\left( p_{k}\right) \in AMIMS$,

$\left( ii\right) $ $p>1$, $0<\alpha <1$, $\left( p_{k}\right) \in AMDMS$
and $\left( n+1\right) =O\left( P_{n}\right) ,$

$\left( iii\right) $ $p>1$, $\alpha =1$ and $\sum\limits_{k=0}^{n-1}\left%
\vert \Delta _{k}\frac{P_{k}}{k+1}\right\vert =O\left( \frac{P_{n}}{n}%
\right) $,

$\left( iv\right) $ $p=1$, $0<\alpha <1$, $\sum\limits_{k=0}^{n-1}\left%
\vert \Delta _{k}p_{k}\right\vert =O\left( n^{-1}\right) $ and $\left(
n+1\right) p_{n}=O\left( P_{n}\right) ,$

$\left( v\right) $ $p=1$, $0<\alpha <1$, $\left( p_{k}\right) \in RBVS$ and $%
\left( n+1\right) =O\left( P_{n}\right) ,$

$\left( vi\right) $ $p=\alpha =1$, $\left( \left( k+1\right) ^{-\beta
}p_{k}\right) \in HBVS$ for some $\beta >0$ and $\left( n+1\right)
p_{n}=O\left( P_{n}\right) $, then%
\begin{equation*}
\left\Vert R_{n}\left( f\right) -f\right\Vert _{L^{p}}=O\left( n^{-\alpha
}\right) .
\end{equation*}
\end{cor}

\begin{rem}
By Theorem 4 we can observe that Theorem 3 and Theorem 1 follow from Remark
1 and Corollary 1 ($\left( i\right) $, $\left( iii\right) $), respectively.
Moreover, since $NDC\subset HBVS$, we can derive from Corollary 1 $\left(
vi\right) $ analogous estimate as in Theorem 2 for the deviation $%
R_{n}\left( f\right) -f$ in the $L^{p}-$norm.
\end{rem}

\section{Auxiliary results}

We shall use the following lemmas for the proof of our theorems:

\begin{lem}
\cite[Theorem 4]{3} If $f\in Lip\left( \alpha ,p\right) $, $p\geq 1$, $%
0<\alpha \leq 1$, then, for any positive integer $n$, $f$ may be
approximated in $L^{p}-$space by a trigonometrical polynomial $t_{n}$ or
order $n$ such that%
\begin{equation*}
\left\Vert f-t_{n}\right\Vert _{L^{p}}=O\left( n^{-\alpha }\right) .
\end{equation*}
\end{lem}

\begin{lem}
\cite[Theorem 5 (i)]{3} If $f\in Lip\left( \alpha ,1\right) $, $0<\alpha <1$%
, then%
\begin{equation*}
\left\Vert \sigma _{n}\left( f\right) -f\right\Vert _{L^{1}}=O\left(
n^{-\alpha }\right) .
\end{equation*}
\end{lem}

\begin{lem}
\cite[p. 541, last line]{3} If $f\in Lip\left( 1,p\right) $ $\left(
p>1\right) $, then%
\begin{equation*}
\left\Vert \sigma _{n}\left( f\right) -S_{n}\left( f\right) \right\Vert
_{L^{p}}=O\left( n^{-1}\right) .
\end{equation*}
\end{lem}

\begin{lem}
\cite[Theorem 6 (i), p 541]{3} Let, for $0<\alpha \leq 1$ and $p>1$, $f\in
Lip\left( \alpha ,p\right) $. Then%
\begin{equation*}
\left\Vert S_{n}\left( f\right) -f\right\Vert _{L^{p}}=O\left( n^{-\alpha
}\right) .
\end{equation*}
\end{lem}

\begin{lem}
Let (\ref{a1}) and (\ref{a2}) hold. If $\left( a_{n,k}\right) \in AMIMS$ or $%
\left( a_{n,k}\right) \in AMDMS$ and $\left( n+1\right) a_{n,0}$, then, for $%
0<\alpha <1$,%
\begin{equation*}
\sum\limits_{k=0}^{n}\left( k+1\right) ^{-\alpha }a_{n,k}=O\left( \left(
n+1\right) ^{-\alpha }\right)
\end{equation*}%
holds.
\end{lem}

\begin{proof}
Let $r=\left[ \frac{n}{2}\right] $. Then, if (\ref{a1}) and (\ref{a2}) hold,%
\begin{equation*}
\sum\limits_{k=0}^{n}\left( k+1\right) ^{-\alpha }a_{n,k}\leq
\sum\limits_{k=0}^{r}\left( k+1\right) ^{-\alpha }a_{n,k}+\left( r+1\right)
^{-\alpha }\sum\limits_{k=r+1}^{n}a_{n,k}
\end{equation*}%
\begin{equation*}
\leq \sum\limits_{k=0}^{r}\left( k+1\right) ^{-\alpha }a_{n,k}+\left(
r+1\right) ^{-\alpha }.
\end{equation*}%
By Abel's transformation, we get%
\begin{equation*}
\sum\limits_{k=0}^{n}\left( k+1\right) ^{-\alpha }a_{n,k}\leq
\sum\limits_{k=0}^{r-1}\left\{ \left( k+1\right) ^{-\alpha }-\left(
k+2\right) ^{-\alpha }\right\} \sum\limits_{i=0}^{k}a_{n,i}
\end{equation*}%
\begin{equation*}
+\left( r+1\right) ^{-\alpha }\sum\limits_{k=0}^{r}a_{n,k}+\left(
r+1\right) ^{-\alpha }\leq \sum\limits_{k=0}^{r-1}\frac{\left( k+2\right)
^{\alpha }-\left( k+1\right) ^{\alpha }}{\left( k+1\right) ^{\alpha
-1}\left( k+2\right) ^{\alpha }}A_{n,k}+\left( r+1\right) ^{-\alpha }.
\end{equation*}%
Using Lagrange's mean value theorem to the function $f\left( x\right)
=x^{\alpha }$ $\left( 0<\alpha <1\right) $ on the interval $\left(
k+1,k+2\right) $ we obtain%
\begin{equation*}
\sum\limits_{k=0}^{n}\left( k+1\right) ^{-\alpha }a_{n,k}\leq
\sum\limits_{k=0}^{r-1}\frac{\alpha }{\left( k+2\right) ^{\alpha }}%
A_{n,k}+\left( r+1\right) ^{-\alpha }.
\end{equation*}%
If $\left( a_{n,k}\right) \in AMIMS$, then%
\begin{equation*}
\sum\limits_{k=0}^{n}\left( k+1\right) ^{-\alpha }a_{n,k}\ll %
{A}_{n,r}\sum\limits_{k=0}^{r-1}\frac{1}{\left( k+2\right) ^{\alpha }}%
+\left( r+1\right) ^{-\alpha }
\end{equation*}%
\begin{equation*}
\ll \left( r+1\right) ^{-\alpha }\sum\limits_{k=0}^{r}a_{n,k}+\left(
r+1\right) ^{-\alpha }\ll \left( n+1\right) ^{-\alpha }.
\end{equation*}%
When $\left( a_{n,k}\right) \in AMDMS$ and $\left( n+1\right)
a_{n,0}=O\left( 1\right) $ we get%
\begin{equation*}
\sum\limits_{k=0}^{n}\left( k+1\right) ^{-\alpha }a_{n,k}\ll %
{A}_{n,0}\sum\limits_{k=0}^{r-1}\frac{1}{\left( k+2\right) ^{\alpha }}%
+\left( r+1\right) ^{-\alpha }
\end{equation*}%
\begin{equation*}
\ll \left( r+1\right) ^{1-\alpha }a_{n,0}+\left( r+1\right) ^{-\alpha }\ll
\left( n+1\right) ^{-\alpha }.
\end{equation*}%
This completes our proof.
\end{proof}

\section{Proofs of the results}

\subsection{Proof of Theorem 4}

$\left( i\right) $ If $\left( a_{n}\right) \in NIS$, then%
\begin{equation*}
\left( n+2\right) \sum\limits_{k=0}^{n}a_{k}=\left( n+1\right)
\sum\limits_{k=0}^{n+1}a_{k}+\sum\limits_{k=0}^{n}a_{k}-\left( n+1\right)
a_{n+1}
\end{equation*}%
\begin{equation*}
\geq \left( n+1\right) \sum\limits_{k=0}^{n+1}a_{k}+\left( n+1\right)
\left( a_{n}-a_{n+1}\right) \geq \left( n+1\right)
\sum\limits_{k=0}^{n+1}a_{k}.
\end{equation*}%
Thus%
\begin{equation*}
\frac{1}{n+2}\sum\limits_{k=0}^{n+1}a_{k}\leq \frac{1}{n+1}%
\sum\limits_{k=0}^{n}a_{k}
\end{equation*}%
and $\left( a_{n}\right) \in NIMS$.

$\left( ii\right) $ Let $\left( a_{n}\right) \in NDS.$ Hence%
\begin{equation*}
\left( n+2\right) \sum\limits_{k=0}^{n}a_{k}=\left( n+1\right)
\sum\limits_{k=0}^{n+1}a_{k}+\sum\limits_{k=0}^{n}a_{k}-\left( n+1\right)
a_{n+1}
\end{equation*}%
\begin{equation*}
\leq \left( n+1\right) \sum\limits_{k=0}^{n+1}a_{k}+\left( n+1\right)
\left( a_{n}-a_{n+1}\right) \leq \left( n+1\right)
\sum\limits_{k=0}^{n+1}a_{k}.
\end{equation*}%
Therefore%
\begin{equation*}
\frac{1}{n+1}\sum\limits_{k=0}^{n}a_{k}\leq \frac{1}{n+2}%
\sum\limits_{k=0}^{n+1}a_{k}
\end{equation*}%
and $\left( a_{n}\right) \in NDMS$.

$\left( iii\right) $ Suppose that $\left( a_{n}\right) \in AMDS$ we have for 
$m\leq l$%
\begin{equation*}
\left( l+1\right) \sum\limits_{i=0}^{m}a_{i}=\left( m+1\right)
\sum\limits_{i=0}^{m}a_{i}+\left( l-m\right) \sum\limits_{i=0}^{m}a_{i}
\end{equation*}%
\begin{equation*}
\geq \left( m+1\right) \left\{ \sum\limits_{i=0}^{m}a_{i}+\frac{1}{K}\left(
l-m\right) a_{m}\right\} \geq \left( m+1\right) \left\{
\sum\limits_{i=0}^{m}a_{i}+\frac{1}{K^{2}}\sum\limits_{i=m+1}^{l}a_{i}%
\right\}
\end{equation*}%
\begin{equation*}
\geq \min \left\{ 1,\frac{1}{K^{2}}\right\} \left( m+1\right)
\sum\limits_{i=0}^{l}a_{i}.
\end{equation*}%
Hence%
\begin{equation*}
\frac{1}{\min \left\{ 1,\frac{1}{K^{2}}\right\} }\frac{1}{m+1}%
\sum\limits_{i=0}^{m}a_{i}\geq \frac{1}{l+1}\sum\limits_{i=0}^{l}a_{i}
\end{equation*}%
and $\left( a_{n}\right) \in AMDMS$.

$\left( iv\right) $ If $\left( a_{n}\right) \in AMIS$, then for $m\leq l$ we
get%
\begin{equation*}
\left( l+1\right) \sum\limits_{i=0}^{m}a_{i}\leq \left( m+1\right) \left\{
\sum\limits_{i=0}^{m}a_{i}+K\left( l-m\right) a_{m}\right\}
\end{equation*}%
\begin{equation*}
\leq \left( m+1\right) \left\{
\sum\limits_{i=0}^{m}a_{i}+K^{2}\sum\limits_{i=m+1}^{l}a_{i}\right\} \leq
\max \left\{ 1.K^{2}\right\} \left( m+1\right) \sum\limits_{i=0}^{l}a_{i}.
\end{equation*}%
Thus%
\begin{equation*}
\frac{1}{m+1}\sum\limits_{i=0}^{m}a_{i}\leq \max \left\{ 1.K^{2}\right\} 
\frac{1}{l+1}\sum\limits_{i=0}^{l}a_{i}
\end{equation*}%
and $\left( a_{n}\right) \in AMIMS$.

The proof is now complete. $\square $

\subsection{Proof of Theorem 5}

We prove the cases $\left( i\right) $ and $\left( ii\right) $ together
utilizing Lemmas 4 and 5. Since%
\begin{equation*}
T_{n}\left( f;x\right) -f\left( x\right)
=\sum\limits_{k=0}^{n}a_{n,k}\left( S_{k}\left( f;x\right) -f\left(
x\right) \right) ,
\end{equation*}%
thus%
\begin{equation*}
\left\Vert T_{n}\left( f\right) -f\right\Vert _{L^{p}}\leq
\sum\limits_{k=0}^{n}a_{n,k}\left\Vert S_{k}\left( f\right) -f\right\Vert
_{L^{p}}\ll \sum\limits_{k=0}^{n}\left( k+1\right) ^{-\alpha
}a_{n,k}=O\left( n^{-\alpha }\right)
\end{equation*}%
and this is (\ref{5}).

Next we consider the case $\left( iii\right) $.

Using two times Abel's transformation and (\ref{a2}) we get that%
\begin{equation*}
T_{n}\left( f;x\right) -f\left( x\right)
=\sum\limits_{k=0}^{n}a_{n,k}\left( S_{k}\left( f;x\right) -f\left(
x\right) \right)
\end{equation*}%
\begin{equation*}
=\sum\limits_{k=0}^{n-1}\left( S_{k}\left( f;x\right) -S_{k+1}\left(
f;x\right) \right) \sum\limits_{i=0}^{k}a_{n,i}+S_{n}\left( f;x\right)
-f\left( x\right)
\end{equation*}%
\begin{equation*}
=S_{n}\left( f;x\right) -f\left( x\right) -\sum\limits_{k=0}^{n-1}\left(
k+1\right) U_{k+1}\left( f;x\right) A_{n,k}
\end{equation*}%
\begin{equation*}
=S_{n}\left( f;x\right) -f\left( x\right) -\sum\limits_{k=0}^{n-2}\left(
A_{n,k}-A_{n,k+1}\right) \sum\limits_{i=0}^{k}\left( i+1\right)
U_{i+1}\left( f;x\right)
\end{equation*}%
\begin{equation*}
-A_{n,n-1}\sum\limits_{k=0}^{n-1}\left( k+1\right) U_{k+1}\left( f;x\right)
=S_{n}\left( f;x\right) -f\left( x\right)
\end{equation*}%
\begin{equation*}
-\sum\limits_{k=0}^{n-2}\left( A_{n,k}-A_{n,k+1}\right)
\sum\limits_{i=0}^{k}\left( i+1\right) U_{i+1}\left( f;x\right) -\frac{1}{n}%
\sum\limits_{i=0}^{n-1}a_{n,i}\sum\limits_{k=0}^{n-1}\left( k+1\right)
U_{k+1}\left( f;x\right) .
\end{equation*}%
Hence%
\begin{equation*}
\left\Vert T_{n}\left( f\right) -f\right\Vert _{L^{p}}\leq \left\Vert
S_{n}\left( f\right) -f\right\Vert _{L^{p}}
\end{equation*}%
\begin{equation}
+\sum\limits_{k=0}^{n-2}\left\vert A_{n,k}-A_{n,k+1}\right\vert \left\Vert
\sum\limits_{i=1}^{k+1}iU_{i}\left( f\right) \right\Vert _{L^{p}}+\frac{1}{n%
}\left\Vert \sum\limits_{k=1}^{n}kU_{k}\left( f;x\right) \right\Vert
_{L^{p}}.  \label{7}
\end{equation}%
Since%
\begin{equation*}
\sigma _{n}\left( f;x\right) -S_{n}\left( f;x\right) =\frac{1}{n+1}%
\sum\limits_{k=1}^{n}kU_{k}\left( f;x\right) ,
\end{equation*}%
thus by Lemma 3%
\begin{equation}
\left\Vert \sum\limits_{k=1}^{n}kU_{k}\left( f\right) \right\Vert
_{L^{p}}=\left( n+1\right) \left\Vert \sigma _{n}\left( f\right)
-S_{n}\left( f\right) \right\Vert _{L^{p}}=O\left( 1\right) .  \label{8}
\end{equation}%
By (\ref{7}), (\ref{8}) and Lemma 4 we get that%
\begin{equation*}
\left\Vert T_{n}\left( f\right) -f\right\Vert _{L^{p}}\ll \frac{1}{n}%
+\sum\limits_{k=0}^{n-1}\left\vert A_{n,k}-A_{n,k+1}\right\vert .
\end{equation*}%
If $\sum\limits_{k=0}^{n-1}\left\vert \Delta _{k}A_{n,k}\right\vert
=O\left( n^{-1}\right) $, then%
\begin{equation*}
\left\Vert T_{n}\left( f\right) -f\right\Vert _{L^{p}}=O\left( n^{-1}\right)
\end{equation*}%
and (\ref{5}) holds.

The cases $\left( iv\right) $ and $\left( v\right) $ we also prove together.

By Abel's transformation%
\begin{equation*}
T_{n}\left( f;x\right) -f\left( x\right)
=\sum\limits_{k=0}^{n}a_{n,k}\left( S_{k}\left( f;x\right) -f\left(
x\right) \right)
\end{equation*}%
\begin{equation*}
=\sum\limits_{k=0}^{n-1}\left( a_{n,k}-a_{n,k+1}\right)
\sum\limits_{i=0}^{k}\left( S_{i}\left( f;x\right) -f\left( x\right)
\right) +a_{n,n}\sum\limits_{k=0}^{n}\left( S_{k}\left( f;x\right) -f\left(
x\right) \right)
\end{equation*}%
\begin{equation*}
=\sum\limits_{k=0}^{n-1}\left( a_{n,k}-a_{n,k+1}\right) \left( k+1\right)
\left( \sigma _{k}\left( f;x\right) -f\left( x\right) \right) +a_{n,n}\left(
n+1\right) \left( \sigma _{n}\left( f;x\right) -f\left( x\right) \right) .
\end{equation*}%
Using Lemma 2 we get%
\begin{equation*}
\left\Vert T_{n}\left( f\right) -f\right\Vert _{L^{1}}\leq
\sum\limits_{k=0}^{n-1}\left\vert a_{n,k}-a_{n,k+1}\right\vert \left(
k+1\right) \left\Vert \sigma _{k}\left( f\right) -f\right\Vert _{L^{1}}
\end{equation*}%
\begin{equation*}
+a_{n,n}\left( n+1\right) \left\Vert \sigma _{n}\left( f\right)
-f\right\Vert _{L^{1}}\ll \sum\limits_{k=0}^{n-1}\left\vert
a_{n,k}-a_{n,k+1}\right\vert \left( k+1\right) ^{1-\alpha }
\end{equation*}%
\begin{equation*}
+a_{n,n}\left( n+1\right) ^{1-\alpha }\leq \left( n+1\right) ^{1-\alpha
}\left( \sum\limits_{k=0}^{n-1}\left\vert a_{n,k}-a_{n,k+1}\right\vert
+a_{n,n}\right) .
\end{equation*}%
When the assumptions $\left( iv\right) $ hold we get%
\begin{equation*}
\left\Vert T_{n}\left( f\right) -f\right\Vert _{L^{1}}=O\left( n^{-\alpha
}\right) .
\end{equation*}%
If $\left( a_{n,k}\right) \in RBVS$, then $\left( a_{n,k}\right) \in AMDS$.
Thus%
\begin{equation*}
\left\Vert T_{n}\left( f\right) -f\right\Vert _{L^{1}}\ll \left( n+1\right)
^{1-\alpha }\left( a_{n,0}+a_{n,n}\right) \ll \left( n+1\right) ^{1-\alpha
}a_{n,0}.
\end{equation*}%
Hence, if $\left( n+1\right) a_{n,0}=O\left( 1\right) $, then (\ref{5})
holds. This ends the proof of the cases $\left( iv\right) $ and $\left(
v\right) $.

Finally, we prove the case $\left( vi\right) $.

Let $t_{n}$ be a trigonometrical polynomial of Lemma 1 of the present paper.
Then for $m\leq n$,%
\begin{equation*}
S_{m}\left( t_{n};x\right) =t_{m}\text{ \ \ \ \ and \ \ \ \ }S_{m}\left(
f;x\right) -t_{m}=S_{m}\left( f-t_{n};x\right) .
\end{equation*}%
Thus%
\begin{equation*}
T_{n}\left( f;x\right) -\sum\limits_{k=0}^{n}a_{n,k}t_{k}\left( x\right)
=\sum\limits_{k=0}^{n}a_{n,k}S_{k}\left( f-t_{n};x\right) ,
\end{equation*}%
where%
\begin{equation*}
S_{k}\left( f-t_{n};x\right) =\frac{1}{\pi }\int\limits_{0}^{2\pi }\left\{
f\left( x+u\right) -t_{n}\left( x+u\right) \right\} \frac{\sin \left( k+%
\frac{1}{2}\right) u}{2\sin \frac{u}{2}}du.
\end{equation*}%
By general form of Minkowski's inequality we get%
\begin{equation*}
\left\Vert T_{n}\left( f\right)
-\sum\limits_{k=0}^{n}a_{n,k}t_{k}\right\Vert _{L^{1}}\leq \frac{1}{2\pi
^{2}}\int\limits_{0}^{2\pi }\left\vert K_{n}\left( u\right) \right\vert
du\int\limits_{0}^{2\pi }\left\vert f\left( x+u\right) -t_{n}\left(
x+u\right) \right\vert dx
\end{equation*}%
\begin{equation*}
=\frac{1}{2\pi ^{2}}\int\limits_{0}^{2\pi }\left\vert K_{n}\left( u\right)
\right\vert \int\limits_{0}^{2\pi }\left\vert f\left( x\right) -t_{n}\left(
x\right) \right\vert dx=\frac{1}{\pi }\left\Vert f-t_{n}\right\Vert
_{L^{1}}\int\limits_{0}^{2\pi }\left\vert K_{n}\left( u\right) \right\vert
du
\end{equation*}%
\begin{equation*}
=\frac{2}{\pi }\left\Vert f-t_{n}\right\Vert _{L^{1}}\int\limits_{0}^{\pi
}\left\vert K_{n}\left( u\right) \right\vert du=\frac{2}{\pi }\left\Vert
f-t_{n}\right\Vert _{L^{1}}\left( \int\limits_{0}^{\pi /n}\left\vert
K_{n}\left( u\right) \right\vert du+\int\limits_{\pi /n}^{\pi }\left\vert
K_{n}\left( u\right) \right\vert du\right)
\end{equation*}%
\begin{equation}
=\frac{2}{\pi }\left\Vert f-t_{n}\right\Vert _{L^{1}}\left(
I_{1}+I_{2}\right) ,  \label{9}
\end{equation}%
where%
\begin{equation*}
K_{n}\left( u\right) =\sum\limits_{k=0}^{n}a_{n,k}\frac{\sin \left( k+\frac{%
1}{2}\right) u}{2\sin \frac{u}{2}}.
\end{equation*}%
Now, we estimate the quantities $I_{1}$ and $I_{2}.$ By (\ref{a2})%
\begin{equation}
I_{1}\ll \int\limits_{0}^{\pi /n}\sum\limits_{k=0}^{n}\left( k+1\right)
a_{n,k}du=O\left( 1\right) .  \label{10}
\end{equation}%
If $\left( \left( k+1\right) ^{-\beta }a_{n,k}\right) \in HBVS$, then $%
\left( \left( k+1\right) ^{-\beta }a_{n,k}\right) \in AMIS$. Hence, for $%
0\leq l\leq m\leq n$,%
\begin{equation*}
Ka_{n,m}\geq a_{n,l}\left( \frac{m+1}{l+1}\right) ^{\beta }\geq a_{n,l}.
\end{equation*}%
Thus $\left( a_{n,k}\right) \in AMIS$. Using this and the assumption $\left(
n+1\right) a_{n,n}=O\left( 1\right) $ we obtain that%
\begin{equation}
I_{2}\ll a_{nn}\int\limits_{\pi /n}^{\pi }u^{-2}du=O\left( 1\right) .
\label{11}
\end{equation}%
Combining (\ref{9})-(\ref{11}) we have%
\begin{equation}
\left\Vert T_{n}\left( f\right)
-\sum\limits_{k=0}^{n}a_{n,k}t_{k}\right\Vert _{L^{1}}\ll \left\Vert
f-t_{n}\right\Vert _{L^{1}}.  \label{12}
\end{equation}%
Further, by using (\ref{12}) and Lemma 1 for $p=\alpha =1$, we get%
\begin{equation*}
\left\Vert T_{n}\left( f\right) -f\right\Vert _{L^{1}}\leq \left\Vert
T_{n}\left( f\right) -\sum\limits_{k=0}^{n}a_{n,k}t_{k}\right\Vert
_{L^{1}}+\left\Vert \sum\limits_{k=0}^{n}a_{n,k}t_{k}-f\right\Vert _{L^{1}}
\end{equation*}%
\begin{equation*}
\ll \frac{1}{n}+\left\Vert \sum\limits_{k=0}^{n}a_{n,k}t_{k}-f\right\Vert
_{L^{1}}\leq \frac{1}{n}+\sum\limits_{k=0}^{n}a_{n,k}\left\Vert
t_{k}-f\right\Vert _{L^{1}}\ll \frac{1}{n}+\sum\limits_{k=0}^{n}\left(
k+1\right) ^{-1}a_{n,k}.
\end{equation*}%
By Abel's transformation%
\begin{equation*}
\left\Vert T_{n}\left( f\right) -f\right\Vert _{L^{1}}\ll \frac{1}{n}%
+\sum\limits_{k=0}^{n-1}\left\vert \frac{a_{n,k}}{\left( k+1\right) ^{\beta
}}-\frac{a_{n,k+1}}{\left( k+2\right) ^{\beta }}\right\vert
\sum\limits_{i=0}^{k}\left( i+1\right) ^{\beta -1}
\end{equation*}%
\begin{equation*}
+\frac{a_{n,n}}{\left( n+1\right) ^{\beta }}\sum\limits_{k=0}^{n}\left(
k+1\right) ^{\beta -1}\ll \frac{1}{n}+\left( n+1\right) ^{\beta
}\sum\limits_{k=0}^{n-1}\left\vert \frac{a_{n,k}}{\left( k+1\right) ^{\beta
}}-\frac{a_{n,k+1}}{\left( k+2\right) ^{\beta }}\right\vert +a_{n,n}.
\end{equation*}%
Since $\left( \left( k+1\right) ^{-\beta }a_{n,k}\right) \in HBVS$ and $%
\left( n+1\right) a_{n,n}=O\left( 1\right) $, then%
\begin{equation*}
\left\Vert T_{n}\left( f\right) -f\right\Vert _{L^{1}}=O\left( n^{-1}\right)
\end{equation*}%
and (\ref{5}) holds.

This completes the proof of Theorem 5. $\square $

\end{document}